\documentclass{article}

\usepackage{arxiv}

\usepackage[utf8]{inputenc} % allow utf-8 input
\usepackage[T1]{fontenc}    % use 8-bit T1 fonts
\usepackage{hyperref}       % hyperlinks
\usepackage{url}            % simple URL typesetting
\usepackage{booktabs}       % professional-quality tables
\usepackage{amsfonts}       % blackboard math symbols
\usepackage{nicefrac}       % compact symbols for 1/2, etc.
\usepackage{microtype}      % microtypography
\usepackage{lipsum}
\usepackage{graphicx}
\graphicspath{ {./images/} }
\usepackage{amsmath}

\newtheorem{theorem}{Théorème}[section]

\newtheorem{definition}[theorem]{Definition}

\newtheorem{proof}[theorem]{Proof}

\newtheorem{lemma}[theorem]{Lemma}

\newtheorem{notation}[theorem]{Notations}

\newtheorem{proposition}[theorem]{Proposition}
\newtheorem{remark}[theorem]{Remark}

\title{Solving a multilevel linear programming problems through a new constructive approach}

\author{
 Mustapha Kaci \\
  Department of Mathematics\\ Faculty of
Mathematics and computer science, Signal Image Parole (SIMPA) Laboratory\\ University of Sciences and Technology of Oran Mohamed Boudiaf USTO-MB, \\El Mnaouar, BP 1505, Bir El Djir 31000, Oran, Algeria. \\
  \texttt{kaci.mustapha.95@gmail.com} \\
   \And
 Sonia Radjef \\
  Department of Mathematics\\ Faculty of
Mathematics and computer science, Signal Image Parole (SIMPA) Laboratory\\ University of Sciences and Technology of Oran Mohamed Boudiaf USTO-MB, \\El Mnaouar, BP 1505, Bir El Djir 31000, Oran, Algeria. \\
  \texttt{soniaradjef@yahoo.fr}  \\
}

\begin{document}
\maketitle
\begin{abstract}
In this paper, an algorithm is developed to solve a  multilevel mono-objective linear programming problem (ML(MO)LPP), where the constructive adaptive method of linear programming is nested. This procedures is the modified versions of the SB. Sinha and S. Sinha's linear programming approach. First, we build a map that reduces the ranges of decision variables that are over control of the previous level's decision maker, depending on the chosen approach, called the range reduction map. Then, we use it to define a new sub-optimality estimate of the adaptive method for the problem under consideration. All the construction stages are carefully checked and illustrated with a numerical example.
\end{abstract}

% keywords can be removed

\section{Introduction}\label{section1}
\par Optimization is a key tool in solving all decision-making problems, whether in science, economics, industry, physics or in other areas. In a classic optimization problem, a decision-maker seeks to select, from among several
decisions, the one that is better or optimal, where the optimality refers to a certain criterion by which
the quality of decisions is measured. However, it has been found that many decision-making problems are
characterized by the presence of multiple decision-makers and/or targets to be optimized \cite{abo,HAN}. Consequently,
traditional optimization approaches cannot be used to formulate and solve these problems because they
cannot produce sufficiently realistic results. New models and methods have therefore been developed to
accommodate these extensions of traditional optimization \cite{kaci,l2}.
\par  Multiobjective optimization techniques have been developed to enable an accurate analysis of trade-offs between competing objectives and to help the decision-maker find an acceptable trade-off \cite{do,hl,ref}. These methods assume that all goals are those of a single decision-maker, which directly impacts their state of well-being. Multiobjective optimization does not take into account the fact that a large number of objectives are categorized into an administrative or hierarchical structure . %\\
\par Multilevel optimization seems to be a very adequate tool for model decision problems where multiple decision makers interact in a hierarchical structure \cite{kaci, abo}.\\
Mathematically, in a multilevel programming problem, the constraints contain a sequence of parametric optimization problems, which must be solved in a predetermined order, and higher level variables are considered parameters in the lower level programming problems.
Multilevel programming  has several successful practical applications in several areas such as supply chain management, network defense, planning, logistics, economics, government, autonomous institutions, agriculture, army, management, schools, hospitals, banks, etc $\ldots$ Although most research on multilevel programming has focused on cases only include two levels (called two-level programming), there are many programming problems that involve more than two levels \cite{l2}.\\
Since the pioneering work of Bracken and McGill \cite{ref36, ref37}, several researchers have published monographs and literature reviews in which theoretical and methodological aspects of two-level optimization were discussed. Several approaches have been developed to solve the multilevel programming problems. For a good bibliography of the multilevel programming problems and their applications, see \cite{kaci,HAN}.  \\
In this work, we consider a ML(MO)LPP (one objective at each level), where the objective  functions and the constraints are linear. We exploit the principle of the adaptive method, developed by R. Gabasov and F.M Kirillova \cite{Gabasov}, to propose a new approach to solve a multilevel linear program.
\par The adaptive method of linear programming is a numerical constructive method, it was described by R. Gabasov, F.M Kirillova and O.I Kostyukova in late 1980s \cite{Gabasov,Gabasov2} to solve the mono-objective linear programming problem  with bounded variables, this method proved its performance  compared to the simplex method, see \cite{Belahcene}. It has been generalized to develop many methods of linear piecewise programming, quadratic programming, quadratic programming \cite{l1}, posynomial programming, optimal control \cite{Gabasov2, Gabasov3} and multiobjective linear programming \cite{delhoum,radjef}.\\
The method is direct, it avoids transforming decision variables and manages the constraints of the problems as they are previously
formulated. This allows to deal with problems in a natural way. So, it generates an important gain in memory space and  permits
 to speed-up the resolution process. The method is really efficient, simple to use and straightforward. In addition, the method incorporates a sub-optimality criterion which makes it possible to stop the algorithm with a desired precision. This could be useful in practical applications.\\
This study is devoted to its application in  multilevel linear  programming problems (ML(L)PP).

%***********************************************************************\\
\par So, in this paper, we propose to improve the approach developed by B. Sinha and S. Sinha in \cite{sinha}, using the adaptive method approach. \\
 First, we construct a map that reduces the range of the decision variables which are controlled by the previous decision-maker, according to the chosen approach called \textit{the range reduction map}, the map is built using the two affine maps given in the Appendix \ref{IRF}. Then we use it to define a new sub-optimality estimate for the problem under consideration. Therefore, we develop two versions of an  algorithm to solve a multilevel linear program. The first version of the algorithm shows how the new sub-optimality estimate occurs, and the second shows how we only use the range reduction map. The algorithms are divided into three parts; in the first part, we   calculate the solution of all the objectives functions independently, then determinate the bounds of all the decision variables, and finally we start the search for the sufficient compromise by solving $P-1$ standard mono-objective linear programming problems with bounded variables using the adaptive method algorithm.
\par The paper is organized as follows. In the next section, we develop the mathematical formulation
of the problem. In section \ref{section2}, a brief description of the adaptive method of linear programming is given. In section \ref{section3}, we construct a map that reduces the ranges of the decision variables according to the multilevel linear programming approach \cite{sinha}. Also, a new sub-optimality estimate of the adaptive method adapted to the multilevel case will be developed. After that, we present a resolution algorithm in section \ref{section4}. Then, the results will be illustrated with  a numerical example in section \ref{section5}.  Finally, a conclusion is given in section \ref{section6}.

\section{Problem formulation}\label{section1}
Consider a $P$-level linear programming problem with single objective function at each level $(P\geq2)$ and  denote DM$_{p}$ the decision maker at $ p^{\text{th}}$ level that has control over the decision variables
 $$
\overline{x}^{p}=x_{p1}, \ldots, x_{pn_{p}}\in\mathbb{R}^{n_{p}},\hspace{0.25cm} p=1,\ldots,P,
$$
 where 
$$
x=(\overline{x}^{1}, \ldots,\overline{x}^{P})\hspace{0.25cm}\text{and}\hspace{0.25cm} n=n_{1}+\ldots+n_{P}.
 $$
 And
\begin{equation}
f_{p}(\overline{x}^{1},\ldots,\overline{x}^{P}):
\begin{array}{ccc}
               \mathbb{R}^{n_{1}}\times\mathbb{R}^{n_{2}}\times \ldots \times\mathbb{R}^{n_{P}}&\longmapsto&
																					                                                                                  \mathbb{R},
																																																														 \end{array}
\end{equation}
where

$$\begin{array}{lll}
f_{p}(x)&=&
c_{p}^{1j}\overline{x}^{1}+c_{p}^{2j}\overline{x}^{2}+\ldots+c_{p}^{Pj}\overline{x}^{P},\hspace{0.3cm} p=\overline{1,P}, j=\overline{1,n_{p}}\vspace{0.2cm}\\
                  &=&c_{p}^{11}x_{11}+\ldots+c_{p}^{1n_{1}}x_{1n_{1}}+c_{p}^{21}x_{21}+\ldots+c_{p}^{2n_{2}}x_{2n_{2}}+c_{p}^{P1}x_{P1}+\ldots+c_{p}^{Pn_{P}}x_{Pn_{P}}
\end{array}
$$
are the coefficients vector of objective  functions of the DM$_{p}$, $p=\overline{1,P}$.\\ The formulation of a $P$-level programming problem with single objective at each level is given as follow:

\begin{equation}\label{mlpp1}
\begin{array}{ll}
 \textbf{Level 1} & \\
& \underset{\overline{x}^{1}}{\max}\hspace{0.1cm} f_{1}(x)=\underset{\overline{x}^{1}}{\max}\hspace{0.1cm}f_{1}(\overline{x}^{1},\ldots,\overline{x}^{P})=\underset{\overline{x}^{1}}{\max}\hspace{0.1cm}c^{t}_{1}x \\
& \text{such that } \overline{x}^{2},\ldots,\overline{x}^{P} \text{solves }\\
 \textbf{Level 2} & \\
& \underset{\overline{x}^{2}}{\max}\hspace{0.1cm} f_{2}(x)=\underset{\overline{x}^{2}}{\max}\hspace{0.1cm}f_{2}(\overline{x}^{1},\ldots,\overline{x}^{P})=\underset{\overline{x}^{2}}{\max}\hspace{0.1cm}c^{t}_{2}x\\
& \vdots\\
& \text{such that } \overline{x}^{P} \text{solves }\\
  \textbf{Level P}  & \\
& \underset{\overline{x}^{P}}{\max}\hspace{0.1cm} f_{P}(x)=\underset{\overline{x}^{P}}{\max}\hspace{0.1cm}f_{P}(\overline{x}^{1},\ldots,\overline{x}^{P})=\underset{\overline{x}^{P}}{\max}\hspace{0.1cm}c^{t}_{P}x
\end{array}
\end{equation}

subject to
$$
x\in S=\left\{x\in\mathbb{R}^{n} : Ax\leq b, x\geq 0, b\in\mathbb{R}^{m}\right\}.
$$
Where $S\neq\emptyset$ is the multilevel convex constraints feasible choice set, $m$ is the number of the constraints, $c_{p}^{ij}$ $i, p=\overline{1,P}$, $j=\overline{1,n_{i}}$ are constants, $b$ is a $m-$vector,  and $A$ is a $m\times n$-matrix.

\section{Preliminaries}\label{section2}
The two basic tools used in this work are the SB. Sinha and S. Sinha's linear approach and the adaptive method of linear programming \cite{Gabasov}. So, in this section, we will briefly present the Adaptive method of linear programming. For the SB. Sinha and S. Sinha's the approach, the reader is returned  to  the reference \cite{sinha}.
\subsection{Adaptive method of linear programming}
Consider the following linear programming problem:
\begin{equation}\label{pl}
    \left\{
    \begin{array}{l}
     \max c^{t} x\\
      Ax = b \\
      l\leq x\leq u
    \end{array},
    \right.
\end{equation}
where $c, x, l, u$ are $n-$vectors, $b$ is $m-$vector, $A=A(I,J)$ is $m\times n-$matrix, $I=\left\{1, \cdots, m\right\}$ and $J=\left\{1, \cdots, n\right\}$. A $n-$vector $x$ satisfying to the main $ Ax= b$ and the direct $  l\leq x\leq u$ constraints of problem (\ref{pl}) is called a feasible solution. The solution $\stackrel{o}{x}$ of problem (\ref{pl}) is called an optimal feasible solution and the feasible solution $x^{\epsilon}$  satisfying the inequality $c^{t} \stackrel{o}{x}-c^{t} x^{\epsilon}\leq \epsilon$ is called $\epsilon-$optimal (or suboptimal) feasible solution. The aim of the algorithm consists in the construction of the suboptimal feasible solution of problem (\ref{pl}) for the given $\epsilon \geq 0$.\\
Its main element is the supporting feasible solution (SFS) $\left\{ x, J_{B}\right\}$, which is a pair from feasible solution $x$ of problem (\ref{pl}) and
the support $J_{B}\subset J$ such that $\left|J_{B}\right|=\left|J\right|=m$ and   $A_{B}=A(I, J_{B})$ being any non-singular submatrix of matrix $A$. A new
feasible solution $\overline{x}$ is constructed in the form $\overline{x}=x+\theta d$ where the direction $d$ is obtained from extremal problem.
The principle of decreasing suboptimality estimation $\beta$ forms the basis of construction of iteration of the adaptive algorithm, it is calculated uniquely
for every supporting feasible solution. As $\beta=\beta(x, J_{B})=\beta(x)+\beta(J_{B})$ where $\beta(x)=c^{t} \stackrel{o}{x}-c^{t} x$ is a measure of feasible
solution nonoptimality, $\beta(J_{B})$ is a measure of support $J_{B}$ nonoptimality, then in contrast with the simplex method, the adaptive method iteration
successfully even in the case $\theta=0$ if $\beta(\overline{J}_{B})< \beta(J_{B})$. The adaptive algorithm is finite at the solution of non-degenerate problems.
\begin{proposition}[Sufficient condition for sub-optimality \cite{Gabasov}]
 Let $\left\{x,J_{B}\right\}$ be a  SFS for the problem $(\ref{pl})$ and $\epsilon$ an arbitrary positive number. If $\beta(x,J_{B})<\epsilon$ then, the feasible solution $x$ is $\epsilon-$optimal.
  \end{proposition}

\section{Linear formulation of the $p^{\text{th}}$-level}\label{section3}
In this section, we will follow the steps (2)-(9) of the linear programming  approach for multilevel programming problem proposed by SB. Sinha and S. Sinha in \cite{sinha}, to formulate the $p^{\text{th}}$-level of the problem  (\ref{mlpp1})  as a linear programming problem. To do, we begin by construct a map $\xi_{p-1}$ that reduce the ranges of the decision variables according to SB. Sinha and S. Sinha's approach, then apply it over the decision variables bounds at each level.
\subsection{Construction of the range reduction map}
\par The higher level DM provides the preferred values of the decision variables under his control to enable the lower level DM to search for his optimum in a narrower feasible region. The basic idea is to reduce the feasible region of a decision variable at each level until a satisfactory point is sought at the last level.

\par The range reduction map $\xi_{p-1}$ that we will construct provides the reduced range of the decision variables that are under the control of the $DM_{p-1}$ exactly as described in the steps (2)-(9) in SB.Sinha and S.Sinha's algorithm, see \cite{sinha}.

\par The idea of the construction is as follows, we use the application 'U' and/or 'L' given in the Appendix \ref{IRF}, to choose the range reduced by the higher level decision-maker, then use the application 'sign' (see Appendix \ref{IRF1}) to choose  between  'U' and 'L'.\\

\begin{notation} $\;$\\

\begin{enumerate}
\item $\stackrel{o}{x}^{p}$ is the DM$_{p}$'s optimal solution, and $\stackrel{o}{x}_{i,j}^{p}$ its $i,j^{\text{th}}$ component.
\item $l^{(p-1)}$, $u^{(p-1)}$ is the lower bound, upper bound of decision variable at the $p^{th}$ level (the new bounds of the decision variables after reducing theirs ranges by the DM$_{(p-1)}$) respectively. And $l^{(p-1)}_{i,j}$, $u^{(p-1)}_{i,j}$ are theirs $i,j^{\text{th}}$ components.
\item The upper right index $^{p}$ indicates that we are in the $p^{\text{th}}$-level, and the lower right indexes $_{i,j}$ refers to the  $i,j^{\text{th}}$ components.
\item 
The functions $L_{p-1,j}$ and $U_{p-1,j}$ are the same to those defined in appendix \ref{IRF}, by taking:
$$
\begin{array}{ccc}
a_{1}=l_{p-1,j}^{(p-1)},& a_{2}=u_{p-1,j}^{(p-1)},& t=\stackrel{o}{x}^{(p-1)}_{p-1,j}.
\end{array}
$$
And the indexes $_{p-1,j}$ means that its acts only on the $p-1,j^{\text{th}}$ components.
\item $LW_{p-1,j}$ and $UP_{p-1,j}$ are the upper range, lower range respectively of the $p-1,j^{\text{th}}$ decision variable’s components.
\end{enumerate}
\end{notation}
\begin{remark}
If $i=p-1$ means that we are dealing with the decisions variable that are under control of the DM$_{p}$.
\end{remark}
\par Consider the $P$-level linear programming problem  (\ref{mlpp1})  where $P\geq2$ with $m$ constraints and $n$ decision variables with one objective function at each level. For $p\in\left\{1,\ldots, P\right\}$, we denote  $\stackrel{o}{x}^{p}$  the optimal solution of the following linear programming problem:
\begin{equation}\label{Pblm1}
\left\{
\begin{array}{l}
     \max \hspace{0.2cm} f_{p}(x)=c_{p}^{t}x           \\
      Ax\leq b\\
x\geq0
\end{array}.
\right.
\end{equation}
Then, we find the bounds of each of the $n$ decision variables as follow
\begin{equation}\label{bound}
l_{i,j}=\underset{q\in\left\{1,\ldots,P\right\}}{\min}\left\{\stackrel{o}{x}_{ij}^{q}\right\}, \, \, \,
u_{i,j}=\underset{q\in\left\{1,\ldots,P\right\}}{\max}\left\{\stackrel{o}{x}_{ij}^{q}\right\}, \, \, \, i=\overline{1,P},j=\overline{1,n_{i}}.
 \end{equation}
We pose $$l^{(1)}_{i,j}=l_{i,j} \, \, \, \text {and} \, \, \,   u^{(1)}_{i,j}=u_{i,j}.$$ Then, we consider the following $2n$ ideals ranges:
\begin{equation}\label{Pblm2}
\begin{array}{lll}
			l^{(1)}_{i,j}\leq x_{i,j} \leq u^{(1)}_{i,j}&\hspace{0.15cm}&i=\overline{1,P},j=\overline{1,n_{i}}.
			\end{array}
\end{equation}

Let $p\in\left\{2, \ldots, P\right\}$ then, the value of each  decision variable  under the control of the $(p-1)^{\text{th}}$-level DM can either be within the
ideal range or at one of its limits, lower or upper.\\ The various cases are:\\

\begin{itemize}
  \item [Case 1:  ] $l^{(p-1)}_{p-1,j}<\stackrel{o}{x}^{(p-1)}_{p-1,j}<u^{(p-1)}_{p-1,j}$.\\
  \item [Case 2:  ] $\stackrel{o}{x}^{(p-1)}_{p-1,j}=u^{(p-1)}_{p-1,j}$.\\
  \item [Case 3:  ] $\stackrel{o}{x}^{(p-1)}_{p-1,j}=l^{(p-1)}_{p-1,j}$.
\end{itemize}

\begin{description}
  \item [\textbf{Case 1}] The $(p-1)^{\text{th}}$  decision maker chooses the reduced range of the decision variables which are under his  control  checking the sign before the decision variable $x_{p-1,j}$, $j=\overline{1,n_{p-1}}$ in $f_{p-1}(x)$.\\ If the sign is negative, we choose  the range
$$
LW_{p-1,j}= \left[l_{p-1,j}^{(p-1)},\stackrel{o}{x}^{(p-1)}_{p-1,j}\right],
$$
using the following map
\begin{equation}
       \begin{array}{llll}
       L_{p-1,j}:& \left[l_{p-1,j}^{(p-1)},u_{p-1,j}^{(p-1)}\right] &  \rightarrow & LW_{p-1,j}            \\
                                                 &  x_{p-1,j}                 &       \rightarrow & \frac{\left(\stackrel{o}{x}^{(p-1)}_{p-1,j}-l_{p-1,j}^{(p-1)}\right)x_{p-1,j}+l_{p-1,j}^{(p-1)}\left(u_{p-1,j}^{(p-1)}-\stackrel{o}{x}^{(p-1)}_{p-1,j}\right)}{u_{p-1,j}^{(p-1)}-l_{p-1,j}^{(p-1)}}
        \end{array}.
\end{equation}
If the sign is positive, we choose  the upper range
$$
UP_{p-1,j} =\left[\stackrel{o}{x}^{(p-1)}_{p-1,j},u_{p-1,j}^{(p-1)}\right],
$$
using the following map:
\begin{equation}
\begin{array}{llll}
U_{p-1,j}:& \left[l_{p-1,j}^{(p-1)},u_{p-1,j}^{(p-1)}\right] & \rightarrow & UP_{p-1,j} \\
                                                 &  x_{p-1,j}                & \rightarrow & \frac{\left(u_{p-1,j}^{(p-1)}-\stackrel{o}{x}^{(p-1)}_{p-1,j}\right)x_{p-1,j}+u_{p-1,j}^{(p-1)}\left(\stackrel{o}{x}^{(p-1)}_{p-1,j}-l_{p-1,j}^{(p-1)}\right)}{u_{p-1,j}^{(p-1)}-l_{p-1,j}^{(p-1)}}
\end{array}.
\end{equation}

\begin{proposition}[Definition]\label{def1}Let $x_{p-1,j}\in\left[l_{p-1,j}^{(p-1)},u_{p-1,j}^{(p-1)}\right]$  and denote
$$
 T_{p-1,j}^{-}=1-sign\left(c^{p-1,j}_{p-1}\right) $$ and $$T_{p-1,j}^{+}=1+sign\left(c^{p-1,j}_{p-1}\right),
$$
 and define
$$
\psi_{p-1,j}\left(x_{p-1,j}\right)=
\frac{T_{p-1,j}^{-}L_{p-1,j}\left(x_{p-1,j}\right)
+T_{p-1,j}^{+}U_{p-1,j}\left(x_{p-1,j}\right)}{2}.
$$

Then, if $sign\left(c^{p-1,j}_{p-1}\right)\neq0$, we have
$$
\psi_{p-1,j}\left(\left[l_{p-1,j}^{(p-1)},u_{p-1,j}^{(p-1)}\right]\right)=\left\{
\begin{array}{lll}
LW_{p-1,j}&\text{if}&c^{p-1,j}_{p-1}<0,\vspace{0.3cm}\\
UP_{p-1,j}&\text{if}&c^{p-1,j}_{p-1}>0.
\end{array}\right.
$$
\end{proposition}

\begin{proof}
Since $sign\left(c^{p-1,j}_{p-1}\right)\neq0$,  then we have $(sign\left(c^{p-1,j}_{p-1}\right)>0\hspace{0.2cm} \text{or}\hspace{0.2cm}  c^{p-1,j}_{p-1}>0)$. So
\begin{itemize}
\item[$\diamond$] If  $sign\left(c^{p-1,j}_{p-1}\right)>0$, then:
$$
\begin{array}{ccc}
T_{p-1,j}^{+}=2&\text{and}&T_{p-1,j}^{-}=0
\end{array}
$$
and
$$
\psi_{p-1,j}\left(\left[l_{p-1,j}^{(p-1)},u_{p-1,j}^{(p-1)}\right]\right)=U_{p-1,j}\left(\left[l_{p-1,j}^{(p-1)},u_{p-1,j}^{(p-1)}\right]\right)=UP_{p-1,j}.
$$

\item[$\diamond$] If  $sign\left(c^{p-1,j}_{p-1}\right)<0$, then 
$$
\begin{array}{ccc}
T_{p-1,j}^{+}=0&\text{and}&T_{p-1,j}^{-}=2
\end{array}
$$
and
$$
\psi_{p-1,j}\left(\left[l_{p-1,j}^{(p-1)},u_{p-1,j}^{(p-1)}\right]\right)=L_{p-1,j}\left(\left[l_{p-1,j}^{(p-1)},u_{p-1,j}^{(p-1)}\right]\right)=LW_{p-1,j}.
$$
\end{itemize}
The case where $sign\left(c^{p-1,j}_{p-1}\right)=0$ is given in proposition \ref{c=0}.
\end{proof}
We can easily see that if $sign\left(c^{p-1,j}_{p-1}\right)\neq0$, then $\psi_{p-1,j}$ sends the ranges $\left[l_{p-1,j}^{(p-1)},u_{p-1,j}^{(p-1)}\right]$ either on $LW_{p-1,j}$ or $UP_{p-1,j}$. Now, if $sign\left(c^{p-1,j}_{p-1}\right)=0$, then the component $x_{p-1,j}$ doesn't appear in $f_{p-1}$, and it doesn't influence the maximization of the considered objective function. So, we won’t make any reduction in the range $\left[l_{p-1,j}^{(p-1)},u_{p-1,j}^{(p-1)}\right]$ (in this case we apply the identity map).\\
 All this settings will be applied only if $c^{p-1,j}_{p}\neq0$. Otherwise, the lower decision-maker (DM$_{p}$) will not use these updates because another time the variable that should receive them does not appear in $f_{p}$. This changes are given by the map $\nu_{p-1,j}$, see  proposition \ref{c=0}.

\begin{proposition}[Definition]\label{c=0}Let $x_{p-1,j}\in\left[l_{p-1,j}^{(p-1)},u_{p-1,j}^{(p-1)}\right]$  and denote
$$
 S_{p-1}^{p-1,j}=sign\left(c^{p-1,j}_{p-1}\right)$$ and $$ S_{p}^{p-1,j}=sign\left(c^{p-1,j}_{p}\right).
$$
Using notation of definition \ref{def1}, we define
$$
\nu_{p-1,j}\left(x_{p-1,j}\right)=\left(S_{p}^{p-1,j}S_{p-1}^{p-1,j}\right)^{2}\psi_{p-1,j}\left(x_{p-1,j}\right)
+\left(1-\left(S_{p}^{p-1,j}S_{p-1}^{p-1,j}\right)^{2}\right)x_{p-1,j}.
$$

Then, if $c^{p-1,j}_{p-1}\neq 0$ and $c^{p-1,j}_{p}\neq0$, we have
$$
\nu_{p-1,j}\left(\left[l_{p-1,j}^{(p-1)},u_{p-1,j}^{(p-1)}\right]\right)=\psi_{p-1,j}\left(\left[l_{p-1,j}^{(p-1)},u_{p-1,j}^{(p-1)}\right]\right).
$$
Otherwise,
$$
\nu_{p-1,j}=Id_{\left[l_{p-1,j}^{(p-1)},u_{p-1,j}^{(p-1)}\right]}.
$$
Where $Id_{\left[l_{p-1,j}^{(p-1)},u_{p-1,j}^{(p-1)}\right]}$ is the identity map.
\end{proposition}
\begin{proof}$\;$\\
\begin{itemize}
\item[$\diamond$]
If $c^{p-1,j}_{p-1}\neq 0$ and $c^{p-1,j}_{p}\neq0$, then
$$
\begin{array}{ccc}
(S_{p}^{p-1,j}S_{p-1}^{p-1,j})^{2}=1 &\text{and}&1-(S_{p}^{p-1,j}S_{p-1}^{p-1,j})^{2}=0
\end{array}
$$
and 
$$
\nu_{p-1,j}\left(\left[l_{p-1,j}^{(p-1)},u_{p-1,j}^{(p-1)}\right]\right)=\psi_{p-1,j}\left(\left[l_{p-1,j}^{(p-1)},u_{p-1,j}^{(p-1)}\right]\right).
$$
\item[$\diamond$] If $c^{p-1,j}_{p-1}=0$ or $c^{p-1,j}_{p}=0$ then
$$
\begin{array}{ccc}
(S_{p}^{p-1,j}S_{p-1}^{p-1,j})^{2}=0&\text{and}&1-(S_{p}^{p-1,j}S_{p-1}^{p-1,j})^{2}=1
\end{array}
$$
and 
$$
\nu_{p-1,j}\left(\left[l_{p-1,j}^{(p-1)},u_{p-1,j}^{(p-1)}\right]\right)=\left[l_{p-1,j}^{(p-1)},u_{p-1,j}^{(p-1)}\right]=Id_{\left[l_{p-1,j}^{(p-1)},u_{p-1,j}^{(p-1)}\right]}.
$$
\end{itemize}
\end{proof}
  \item [\textbf{Case 2}] If $\stackrel{o}{x}^{(p-1)}_{p-1,j}=u^{(p-1)}_{p-1,j}$, then the choice of reduction will be as follow
$$
LW_{p-1,j}^{\alpha_{p-1,j}}=\left[l_{p-1,j}^{(p-1)},u_{p-1,j}^{(p-1)}-\alpha_{p-1,j}\right],
$$
using the map
\begin{equation}
\begin{array}{llll}
L_{p-1,j}^{\alpha_{p-1,j}}:& \left[l_{p-1,j}^{(p-1)},u_{p-1,j}^{(p-1)}\right] & \rightarrow & LW_{p-1,j}^{\alpha_{p-1,j}}\\
                                                 &  x_{p-1,j}                         & \rightarrow & \frac{\left(u_{p-1,j}^{(p-1)}-\left(l_{p-1,j}^{(p-1)}+\alpha_{p-1,j}\right)\right)x_{p-1,j}+\alpha_{p-1,j} l_{p-1,j}^{(p-1)}}{u_{p-1,j}^{(p-1)}-l_{p-1,j}^{(p-1)}}
\end{array}.
\end{equation}
\item [\textbf{Case 3}]   If $\stackrel{o}{x}^{p-1,j}_{p-1}=l^{(p-1)}_{p-1,j}$  then, we choose the range

$$
UP_{p-1,j}^{\alpha_{p-1,j}}=\left[l_{p-1,j}^{(p-1)}+\alpha_{p-1,j},u_{p-1,j}^{(p-1)}\right],
$$
using the map
\begin{equation}
\begin{array}{llll}
U_{p-1,j}^{\alpha_{p-1,j}}:& \left[l_{p-1,j}^{(p-1)},u_{p-1,j}^{(p-1)}\right] & \rightarrow & UP_{p-1,j}^{\alpha_{p-1,j}}\\
                                                 &  x_{p-1,j}                         & \rightarrow & \frac{\left(u_{p-1,j}^{(p-1)}-\left(l_{p-1,j}^{(p-1)}+\alpha_{p-1,j}\right)\right)x_{p-1,j}+\alpha_{p-1,j}u_{p-1,j}^{(p-1)}}{u_{p-1,j}^{(p-1)}-l_{p-1,j}^{(p-1)}}
\end{array}.
\end{equation}
Where  $\alpha_{p-1,j}\in\mathbb{R}$ is an arbitrary number chosen by the DM$_{p-1}$, such that
\begin{equation}\label{alpha}
\alpha_{p-1,j}<u_{p-1,j}^{(p-1)}-l_{p-1,j}^{(p-1)}.
\end{equation}
\begin{proposition}[Definition]\label{def2}
 Let's denote
$$A_{p-1,j}=sign\left(u^{(p-1)}_{p-1,j}-\stackrel{o}{x}^{(p-1)}_{p-1,j}\right)$$ and $$B_{p-1,j}=sign\left(\stackrel{o}{x}^{(p-1)}_{p-1,j}-l^{(p-1)}_{p-1,j}\right).
$$
And define
$$
\hat{\psi}_{p-1,j}(x_{p-1,j})=B_{p-1,j}L_{p-1,j}^{\alpha_{p-1,j}}(x_{p-1,j})
+A_{p-1,j}U_{p-1,j}^{\alpha_{p-1,j}}(x_{p-1,j}).
$$ 
Then, for all $x_{p-1,j}\in\left[l_{p-1,j}^{(p-1)},u_{p-1,j}^{(p-1)}\right]$, we have
\begin{itemize}
\item[$\diamond$] If $sign\left(u^{(p-1)}_{p-1,j}-\stackrel{o}{x}^{(p-1)}_{p-1,j}\right)=0$, then
$$
\hat{\psi}_{p-1,j}\left(x_{p-1,j}\right)=L_{p-1,j}^{\alpha_{p-1,j}}(x_{p-1,j}).
$$

\item[$\diamond$]If $sign\left(\stackrel{o}{x}^{(p-1)}_{p-1,j}-l^{(p-1)}_{p-1,j}\right)=0$, then
$$
\hat{\psi}_{p-1,j}\left(x_{p-1,j}\right)=U_{p-1,j}^{\alpha_{p-1,j}}(x_{p-1,j}).
$$
\end{itemize}

\end{proposition}
\begin{proof}Immediate.
\end{proof}

\begin{proposition}[Definition]\label{c=1}Let $x_{p-1,j}\in\left[l_{p-1,j}^{(p-1)},u_{p-1,j}^{(p-1)}\right]$,
using notation of proposition \ref{c=0}, we define
$$
\hat{\nu}_{p-1,j}\left(x_{p-1,j}\right)=\left(S_{p}^{p-1,j}S_{p-1}^{p-1,j}\right)^{2}\hat{\psi}_{p-1,j}\left(x_{p-1,j}\right)
+\left(1-\left(S_{p}^{p-1,j}S_{p-1}^{p-1,j}\right)^{2}\right)x_{p-1,j}.
$$

Then, if $c^{p-1,j}_{p-1}\neq 0$ and $c^{p-1,j}_{p}\neq0$, we have
$$
\hat{\nu}_{p-1,j}\left(\left[l_{p-1,j}^{(p-1)},u_{p-1,j}^{(p-1)}\right]\right)=\hat{\psi}_{p-1,j}\left(\left[l_{p-1,j}^{(p-1)},u_{p-1,j}^{(p-1)}\right]\right).
$$
Otherwise,
$$
\hat{\nu}_{p-1,j}=Id_{\left[l_{p-1,j}^{(p-1)},u_{p-1,j}^{(p-1)}\right]}.
$$
Where $Id_{\left[l_{p-1,j}^{(p-1)},u_{p-1,j}^{(p-1)}\right]}$ is the identity map.
\end{proposition}
\begin{proof}
The proof is similar to that of the proposition \ref{c=0}.
\end{proof}

\end{description}
\begin{notation}
Using the notations of proposition \ref{def2}, let denote
$$
\begin{array}{l}
A^{+}_{p-1,j}=1+A_{p-1,j}, A^{-}_{p-1,j}=1-A_{p-1,j},\vspace{0.25cm}\\
B^{+}_{p-1,j}=1+B_{p-1,j}, B^{-}_{p-1,j}=1-B_{p-1,j}
\end{array}
$$
and
$$
\begin{array}{l}
H_{p-1,j}^{1}=A^{+}_{p-1,j}A^{-}_{p-1,j}+B^{+}_{p-1,j}B^{-}_{p-1,j}, \vspace{0.25cm}\\
 H_{p-1,j}^{2}=-A_{p-1,j}B_{p-1,j}. \vspace{0.25cm}
\end{array}
$$
\end{notation}

\begin{proposition}[Range reduction map]\label{RRM}

Let $i=\overline{2,P}$, $j=\overline{1,n_{i}}$ and $ x_{ij}\in\left[l^{(p-1)}_{i,j},u^{(p-1)}_{i,j}\right]$, and define
 \begin{equation}
\xi_{p-1}(x_{ij})=   \left\{\begin{array}{lll}
                              H_{p-1,j}^{1}\hat{\psi}_{i,j}(x_{ij})+H_{p-1,j}^{2}\hat{\nu}_{i,j}(x_{ij})&\text{if}&i=p-1 \\
                               x_{ij}                                                                                    &\text{if} & i\neq p-1
                        \end{array}
                        \right.
 \end{equation}
and
$$
\xi_{0}\left(x_{ij}\right)=x_{ij}.
$$
Then, the  map $\xi_{p-1}$ will reduce the ranges of the decision variable $x_{ij}$ with regard to all the previous settings.
\end{proposition}

\begin{proof}Let $i=\overline{2,P}$, $j=\overline{1,n_{i}}$ and $ x_{ij}\in\left[l^{(p-1)}_{i,j},u^{(p-1)}_{i,j}\right]$, then
\begin{description}
\item[Case 1]If $l^{(p-1)}_{p-1,j}<\stackrel{o}{x}^{(p-1)}_{p-1,j}<u^{(p-1)}_{p-1,j}$, we have
$$
\begin{array}{lclcl}
 A_{p-1,j}=sign\left(u^{(p-1)}_{p-1,j}-\stackrel{o}{x}^{(p-1)}_{p-1,j}\right)=1&\Rightarrow&A^{+}_{p-1,j}=1+A_{p-1,j}=2&\text{and}&A^{-}_{p-1,j}=1-A_{p-1,j}=0,\\\vspace{0.2cm}
B_{p-1,j}=sign\left(\stackrel{o}{x}^{(p-1)}_{p-1,j}-l^{(p-1)}_{p-1,j}\right)=-1&\Rightarrow&B^{+}_{p-1,j}=1+B_{p-1,j}=0&\text{and}&B^{-}_{p-1,j}=1-B_{p-1,j}=2.
\end{array}
$$
And
$$
\begin{array}{lclcl}
H_{p-1,j}^{1}=A^{+}_{p-1,j}A^{-}_{p-1,j}+B^{+}_{p-1,j}B^{-}_{p-1,j}&=&2\times 0+0\times 2&=&0, \vspace{0.25cm}\\
 H_{p-1,j}^{2}=-A_{p-1,j}B_{p-1,j}&=&-1\times -1&=&1.
\end{array}
$$
Then, we get
$$
\xi_{p-1}(x_{ij})=   \left\{\begin{array}{lclcl}
                              0\times\hat{\psi}_{i,j}(x_{ij})+\hat{\nu}_{i,j}(x_{ij})&=&\hat{\nu}_{i,j}(x_{ij})&\text{if}&i=p-1, \\
                               x_{ij}         &&                                                                           &\text{if} & i\neq p-1.
                        \end{array}
                        \right.
$$

  \item [Case 2] $\stackrel{o}{x}^{(p-1)}_{p-1,j}=u^{(p-1)}_{p-1,j}$, we have
$$
\begin{array}{lclcl}
 A_{p-1,j}=sign\left(u^{(p-1)}_{p-1,j}-\stackrel{o}{x}^{(p-1)}_{p-1,j}\right)=0&\Rightarrow&A^{+}_{p-1,j}=1+A_{p-1,j}=1&\text{and}&A^{-}_{p-1,j}=1-A_{p-1,j}=1,\\\vspace{0.2cm}
B_{p-1,j}=sign\left(\stackrel{o}{x}^{(p-1)}_{p-1,j}-l^{(p-1)}_{p-1,j}\right)=1&\Rightarrow&B^{+}_{p-1,j}=1+B_{p-1,j}=2&\text{and}&B^{-}_{p-1,j}=1-B_{p-1,j}=0.
\end{array}
$$
And
$$
\begin{array}{lclcl}
H_{p-1,j}^{1}=A^{+}_{p-1,j}A^{-}_{p-1,j}+B^{+}_{p-1,j}B^{-}_{p-1,j}&=&1\times 1+2\times 0&=&1, \vspace{0.25cm}\\
 H_{p-1,j}^{2}=-A_{p-1,j}B_{p-1,j}&=&-0\times 1&=&0.
\end{array}
$$
Then, we get
$$
\xi_{p-1}(x_{ij})=   \left\{\begin{array}{lclcl}
                              \hat{\psi}_{i,j}(x_{ij})+0\times\hat{\nu}_{i,j}(x_{ij})&=&\hat{\psi}_{i,j}(x_{ij})&\text{if}&i=p-1, \\
                               x_{ij}         &&                                                                           &\text{if} & i\neq p-1.
                        \end{array}
                        \right.
$$

  \item [Case 3] If $\stackrel{o}{x}^{(p-1)}_{p-1,j}=l^{(p-1)}_{p-1,j}$, we have
$$
\begin{array}{lclcl}
 A_{p-1,j}=sign\left(u^{(p-1)}_{p-1,j}-\stackrel{o}{x}^{(p-1)}_{p-1,j}\right)=1&\Rightarrow&A^{+}_{p-1,j}=1+A_{p-1,j}=2&\text{and}&A^{-}_{p-1,j}=1-A_{p-1,j}=0,\\\vspace{0.2cm}
B_{p-1,j}=sign\left(\stackrel{o}{x}^{(p-1)}_{p-1,j}-l^{(p-1)}_{p-1,j}\right)=0&\Rightarrow&B^{+}_{p-1,j}=1+B_{p-1,j}=1&\text{and}&B^{-}_{p-1,j}=1-B_{p-1,j}=1.
\end{array}
$$
And
$$
\begin{array}{lclcl}
H_{p-1,j}^{1}=A^{+}_{p-1,j}A^{-}_{p-1,j}+B^{+}_{p-1,j}B^{-}_{p-1,j}&=&2\times 0+1\times 1&=&1, \vspace{0.25cm}\\
 H_{p-1,j}^{2}=-A_{p-1,j}B_{p-1,j}&=&-1\times 0&=&0.
\end{array}
$$
Then, we get
$$
\xi_{p-1}(x_{ij})=   \left\{\begin{array}{lclcl}
                              \hat{\psi}_{i,j}(x_{ij})+0\times\hat{\nu}_{i,j}(x_{ij})&=&\hat{\psi}_{i,j}(x_{ij})&\text{if}&i=p-1, \\
                               x_{ij}         &&                                                                           &\text{if} & i\neq p-1.
                        \end{array}
                        \right.
$$

\end{description}

\end{proof}
\begin{remark}
Let $p=2,\ldots P$, then the range reduction map $\xi_{p-1}$ reduce the ranges of only components $x_{p-1,j}$ for all $j=1,\ldots,n_{p}$.
\end{remark}

\begin{definition}
The linear formulation of the $p^{\text{th}}$-level is given as follow:
	\begin{equation}\label{model}
	\left\{
\begin{array}{l}
	            \underset{x}{\max} f_{p}(x)=c^{t}_{p}x\\
       \begin{array}{l}
                                             Ax\leq b                 \\
		                                         \xi_{p-1}\left(l^{(p-1)}\right)\leq x\leq \xi_{p-1}\left(u^{(p-1)}\right)\\
		                                         x\geq 0
		                             \end{array}
\end{array},\right.
\end{equation}
where
\begin{equation}\label{equation1}
\xi_{p-1}(l^{(p-1)})=\left(\xi_{p-1}\left(l^{(p-1)}_{ij}\right), i=\overline{1,P},j=\overline{1,n_{i}}\right)
\end{equation}
and
\begin{equation}\label{equation2}
\xi_{p-1}(u^{(p-1)})=\left(\xi_{p-1}\left(u^{(p-1)}_{ij}\right), i=\overline{1,P},j=\overline{1,n_{i}}\right)
\end{equation}.
\end{definition}
\subsection{Sub-optimality estimate for the  $p^{\text{th}}$-level linear programming problem}
In order to pass from an indexing by two indexes to an indexing by one index, we give the following definition
\begin{definition}\textit{Index change map}\\
Let $I=\left\{1, \ldots, n\right\}$, $I_{1}=\left\{1, \ldots, n_{1}\right\}$  and for all $i=\overline{2,P}$, we denote
$$
 I_{i}=\left\{ \sum_{t=1}^{i-1}{n_{t}}+j\hspace{0.2cm}/ \hspace{0.2cm}j=\overline{1,n_{i}}\right\}.
$$
And for all $i \in\left\{1, \ldots, P\right\}$, we define
$$
\begin{array}{llll}
e:& I_{i}&\rightarrow&I\\
    & j&\rightarrow&e_{i,j}:=e(i,j)=\left\{\begin{array}{lll}
		                              \sum_{t=1}^{i-1}n_{t} +j&\text{if}&i>1\\
																    j                                        &\text{if}&i=1
												 \end{array} \right.
\end{array}.
$$

\end{definition}

\begin{definition}
Let $\left\{x^{p},J_{B}^{p}\right\}$ be a  SFS for the problem $(\ref{model})$, and define:
\begin{equation}\label{formulee0}
\begin{array}{lll}
\beta_{p}^{\xi}(x^{p},J_{B}^{p})=&
\sum_{E_{e_{i,j}}^{p}>0,e_{i,j}\in J_{N}^{p}} E_{e_{i,j}}^{p}\left(\xi_{p-1}\left(x^{p}_{e_{i,j}}\right)-\xi_{p-1}\left(l^{(p-1)}_{e_{i,j}}\right)\right)\\
&  +
\sum_{E_{e_{i,j}}^{p}<0,e_{i,j}\in J_{N}^{p}} E_{e_{i,j}}^{p}\left(\xi_{p-1}\left(x^{p}_{e_{i,j}}\right)-\xi_{p-1}\left(u^{(p-1)}_{e_{i,j}}\right)\right).

\end{array}
\end{equation}
$\beta_{p}^{\xi}(x^{p},J_{B}^{p})$ is the sub-optimality estimate for the $p^{\text{th}}$-level.
\end{definition}

\begin{proposition}[Sufficient condition for sub-optimality of the $p^{\text{th}}$-level]
 Let $\left\{x^{p},J_{B}^{p}\right\}$ be a  SFS for the problem $(\ref{model})$ and $\epsilon_{p}$ an arbitrary positive number. If $\beta_{p}^{\xi}(x^{p},J_{B}^{p})<\epsilon_{p}$ then, the feasible solution $x^{p}$ is $\epsilon_{p}-$optimal.
  \end{proposition}
	\begin{proof} See \cite{Gabasov}.
	\end{proof}
	\begin{definition}
	The optimal solution $\stackrel{c}{x}^{P}$  of the linear programming problem $(\ref{model})$ is called  optimal satisfactory compromise of the multilevel linear program $(\ref{mlpp1})$.
	\end{definition}
	\begin{proposition}[A sufficient condition for satisfactory compromise of a multilevel linear program]
 Let $\left\{\stackrel{c}{x}^{P},J_{B}^{P}\right\}$ be a  SFS for the problem $(\ref{model})$ and $\epsilon_{P}$ an arbitrary positive number. If $\beta_{P}^{\xi}(\stackrel{c}{x}^{P},J_{B}^{P})<\epsilon_{P}$ then, the feasible solution $\stackrel{c}{x}^{P}$ is an $\epsilon_{P}-$optimal satisfactory compromise of the multilevel linear program $(\ref{mlpp1})$.
  \end{proposition}
	\begin{proof} The satisfactory compromise of the problem $(\ref{mlpp1})$ is obtained as a solution of the linear programming problem with bounded variable (\ref{model}) using the adaptive method at the  iteration number $P-1$, then See \cite{Gabasov}.
	\end{proof}

	\section{Adaptive method algorithm for a ML(MO)LPP}\label{section4}
	\par Using the above construction, we give an algorithm to solve the multilevel linear program $(\ref{mlpp1})$. The algorithm follow  the scheme:
	\begin{itemize}
	\item Calculate the solution of all the objectives independently.
	\item Determinate the bounds of all the decision variables.
	\item Start the search for the sufficient/compromise by solving $P-1$ standard linear programming problems with bounded variables using the adaptive method algorithm.
	\end{itemize}

\subsection{Algorithm}	
Consider the following linear programming problem:
	\begin{equation}\label{model2}
	\left\{
\begin{array}{l}
	            \underset{x}{\max} f_{p}(x)=c^{t}_{p}x\\
              \begin{array}{l}
                                             Ax\leq b                 \\
		                                         l^{(p)}\leq x\leq u^{(p)}\\
		                                         x\geq 0
		                             \end{array}
\end{array},\right.
\end{equation}
Then, the adaptive method to solve the ML(MO)LPP $(\ref{mlpp1})$ is given by the following steps
	\begin{description}
	\item[Step 1] Set $p=1.$
	\item[Step 2] Explicit the model (\ref{Pblm1}).
	\item[Step 3] Get the solutions $\stackrel{o}{x}^{p}$ of the model (\ref{Pblm1}) "solution of the $p^{\text{th}}$-level".
	\item[Step 4] $p=p+1.$
	\item[Step 5] If $p>P$  then stop with $P$ optimal solutions $\stackrel{o}{x}^{1}$,\ldots,$\stackrel{o}{x}^{P}$ and go to \textbf{Step 6}. \\Otherwise, go to \textbf{Step 2}.
		\item[Step 6 ] Get the bounds of all the decision  variables $l^{(1)}_{e_{ij}}$  and $u^{(1)}_{e_{ij}}$, $i\in\left\{1,P\right\}, j\in\left\{1,n_{i}\right\}, e_{ij}\in\left\{1,n\right\}$ using the formulas (\ref{bound}).
		\item[Step 7] Set $p=2.$
		\item[Step 8] Compute $l^{(p)}=\xi_{p-1}\left(l^{(p-1)}\right)$ and $u^{(p)}=\xi_{p-1}\left(u^{(p-1)}\right)$.
		\item[Step 9] Explicit the model (\ref{model2}).
		\item[Step 10]  Solve the model (\ref{model2}), using the adaptive method algorithm given in section \ref{section2}, to get the optimal solution $x^{p}$.
		\item[Step 11] Set $\stackrel{c}{x}^{p}=x^{p}$.
		 \item[Step 12] $p=p+1.$
		\item[Step 13] If $p>P$  then, stop with a satisfactory solution $\stackrel{c}{x}^{p}$  to the multilevel programming problem (\ref{mlpp1}). \\ Otherwise, go to \textbf{Step 8}.
		\end{description}
	\subsubsection{Algorithm explanation}
	\begin{itemize}
\item[$\diamond$] \textbf{Step 1-Step 5}: Calculate the optimal solution of all the objective functions that constitutes the $P-$levels problem \ref{mlpp1} independently.\vspace{0.2cm}
\item[$\diamond$] \textbf{Step 6}: Illustrate the initial bounds of the decision variables vector.\vspace{0.2cm}
\item[$\diamond$] \textbf{Step 7-Step 13}: The algorithm stop with a satisfactory solution $\stackrel{c}{x}^{P}$ after  $P-1$ iterations, So in each iteration the $DM_{p}$ chose the reduced range for the decision variables that are under his control using the \textit{range reduction map} \ref{RRM}, then solve its mono-objective linear programming using the adaptive method described in section \ref{section2} and  \cite{Gabasov} .
\item[$\diamond$] We apply the range reduction map $\xi_{p-1}$ only on the bounds of the decision variables $l^{(p-1)}$ and $u^{(p-1)}$ and apply the adaptive method algorithm with the old sub-optimality estimate.
	\end{itemize}
	
	\section{Numerical example }\label{section5}
	
Consider the following three-level mono-objective linear programming problem

 \begin{equation}\label{exemple131}
\begin{array}{ll}
  \textbf{Level 1} &\\
& \underset{x_{1},x_{2}}{\max}\hspace{0.1cm} f_{1}(x)=-5x_{1}+x_{2}+2x_{3}+3x_{4} \\
&  x_{3},x_{4} \hspace{0.1cm} \text{solves }\\
 \textbf{Level 2} & \\
 & \underset{x_{3}}{\max}\hspace{0.1cm} f_{2}(x)=6x_{1}+2x_{2}-3x_{3}+x_{4}\\
&  x_{4} \hspace{0.1cm} \text{solves }\\
 \textbf{Level 3}&\\
& \underset{x_{4}}{\max}\hspace{0.1cm} f_{3}(x)=-x_{2}+2x_{3}+3x_{4},
\end{array}
\end{equation}

subject to
$$
x\in S=\left\{x\in\mathbb{R}^{n} / Ax\leq b, x\geq 0, b\in\mathbb{R}^{m}\right\}\neq\emptyset,
$$
where
$$
\begin{array}{cc}
A=\left(\begin{array}{rrrr}
3&2&1&2\\
1&3&-2&-1\\
-2&-5&-2&0\\
-1&4&1&0\\
1&-1&1&2\\
1&0&1&3\\
0&0&0&1
\end{array}\right),&\hspace{0.3cm}
b=\left(\begin{array}{r}6\\ 3\\ -2\\ 2 \\5\\ 4 \\2
\end{array}\right)
\end{array}.
$$
We get the optimal solution: $$
\begin{array}{lcl}
  f_{1}=6  &  \text{at \, \, }& \stackrel{o}{x}^{1}=\left(\begin{array}{cccc}0&0&2&\frac{2}{3}
  \end{array}\right),\vspace{0.25cm}\\
  f_{2}=12 & \text{at \, \, }& \stackrel{o}{x}^{2}=\left(\begin{array}{cccc}2&0&0&0 \end{array}\right),\vspace{0.25cm}\\
  f_{3}=6 & \text{at \, \, }&\stackrel{o}{x}^{3}=\left(\begin{array}{cccc}1&0&3&0 \end{array}\right).
\end{array}
$$ Then, the bounds of the decision  variables are
	$$
	\begin{array}{llll}
	l^{(1)}_{11}=0,&l^{(1)}_{12}=0,&l^{(1)}_{21}=0,&l^{(1)}_{31}=0,\vspace{0.25cm}\\
	u^{(1)}_{11}=2,&u^{(1)}_{12}=0,&u^{(1)}_{21}=3,&u^{(1)}_{31}=\frac{2}{3}.
	\end{array}
	$$
	\begin{description}
	\item[Iteration 1]
	 By taking $\alpha_{11}=0.25$, $\alpha_{12}=\alpha_{21}=\alpha_{31}=0$  and
	$$
	\begin{array}{llll}
	 \xi_{1}\left(l^{(1)}_{11}\right)=U_{11}^{0.25}\left(l^{(1)}_{11}\right)&=0.25,&\xi_{1}\left(u^{(1)}_{11}\right)=Id_{[0,2]}\left(u^{(1)}_{11}\right)&=2,\vspace{0.25cm}\\
	 \xi_{1}\left(l^{(1)}_{12}\right)=Id_{\left\{0\right\}}\left(l^{(1)}_{12}\right)&=0,&\xi_{1}\left(u^{(1)}_{12}\right)=Id_{\left\{0\right\}}\left(u^{(1)}_{12}\right)&=0,\vspace{0.25cm}\\
	\xi_{1}\left(l^{(1)}_{21}\right)=Id_{[0,3]}\left(l^{(1)}_{21}\right)&=0,&\xi_{1}\left(u^{(1)}_{21}\right)=Id_{[0,3]}\left(u^{(1)}_{21}\right)&=3,\vspace{0.25cm}\\
	 \xi_{1}\left(l^{(1)}_{31}\right)=Id_{\left[0,\frac{2}{3}\right]}\left(l^{(1)}_{31}\right)&=0&\xi_{1}\left(u^{(1)}_{31}\right)=Id_{\left[0,\frac{2}{3}\right]}\left(u^{(1)}_{31}\right)&=\frac{2}{3}.
	
	\end{array}
	$$
We put  
$$
\begin{array}{lcl}
\begin{array}{lcl}
d^{(2)-}&=&\left(\begin{array}{cccc}
                        0.25&0&0&0
                       \end{array}
\right),
\end{array}&\text{and}&\begin{array}{lcl}
d^{(2)+}&=&\left(\begin{array}{cccc}2&0&3&\frac{2}{3}\end{array}\right).
\end{array}
\end{array}
$$ 

 Then, we consider the following linear programming problem:
		\begin{equation}\label{sol1}
	\left\{
\begin{array}{l}
	            \underset{x}{\max} f_{2}(x)=6x_{1}+2x_{2}-3x_{3}+x_{4}\\
                                             Ax\leq b                 \\
		                                         l^{(2)}\leq x\leq u^{(2)}\\
		                                         x\geq 0
		                             \end{array}.\right.
\end{equation}
Using the adaptive method algorithm , we get the optimal solution of the linear problem (\ref{sol1})
$$
\stackrel{c}{x}^{2}=\left(\begin{array}{cccc}2&0&0&0\end{array}\right).
$$
\item[Iteration 2] We set
	$$\alpha_{21}=0.25, \, \, \, \alpha_{11}=\alpha_{12}=\alpha_{31}=0$$  and
	$$
	\begin{array}{llll}
	 l^{(3)}_{11}=\xi_{2}\left(l^{(2)}_{11}\right)=Id_{[0.25,2]}\left(l^{(2)}_{11}\right)&=0.25,&u^{(3)}_{11}=\xi_{2}\left(u^{(2)}_{11}\right)=Id_{[0.25,2]}\left(u^{(2)}_{11}\right)&=2,\vspace{0.25cm}\\
	 l^{(3)}_{12}=\xi_{2}\left(l^{(2)}_{12}\right)=Id_{\left\{0\right\}}\left(l^{(2)}_{12}\right)&=0,&u^{(3)}_{12}=\xi_{2}\left(u^{(2)}_{12}\right)=Id_{\left\{0\right\}}\left(u^{(2)}_{12}\right)&=0,\vspace{0.25cm}\\
	 l^{(3)}_{21}=\xi_{2}\left(l^{(2)}_{21}\right)=U^{0.25}_{21}\left(l^{(2)}_{21}\right)&=0.25,&u^{(3)}_{21}=\xi_{2}\left(u^{(2)}_{21}\right)=Id_{[0,3]}\left(u^{(2)}_{21}\right)&=3,\vspace{0.25cm}\\
	 l^{(3)}_{31}=\xi_{2}\left(l^{(2)}_{31}\right)=Id_{\left[0,\frac{2}{3}\right]}\left(l^{(2)}_{31}\right)&=0&u^{(3)}_{31}=\xi_{2}\left(u^{(2)}_{31}\right)=Id_{\left[0,\frac{2}{3}\right]}\left(u^{(2)}_{31}\right)&=\frac{2}{3}.
	
	\end{array}
	$$
	Then, we put  
	$$
	\begin{array}{lcl}
	l^{(3)}=\left(\begin{array}{cccc}0.25&0&0.25&0\end{array}\right)&
	\text{and} &
	u^{(3)}=\left(\begin{array}{cccc}2&0&3&\frac{2}{3}\end{array}\right),
	\end{array}
	$$
	and we consider the following linear programming problem:
		\begin{equation}\label{sol2}
	\left\{
\begin{array}{l}%\label{sol2}
	            \underset{x}{\max} f_{3}(x)=-x_{2}+2x_{3}+3x_{4}\\
                                             Ax\leq b                 \\
		                                         l^{(3)}\leq x\leq u^{(3)}\\
		                                         x\geq 0
		                    \end{array}.\right.
\end{equation}
	Finally, using the adaptive method algorithm of mono-objective linear programming, we get the satisfactory/compromise of the three-level linear programming problem (\ref{exemple131}) as a solution of the problem (\ref{sol2}), equal to
		$$
		\stackrel{c}{x}^{3}=\left(\begin{array}{cccc}0.2500&0&2.7671&0\end{array}\right).
		$$
		
	Finally, we get
		$$
		\begin{array}{ccc}
		f_{1}\left(\stackrel{c}{x}^{3}\right)=4.2841,&f_{2}\left(\stackrel{c}{x}^{3}\right)=-6.8012,&f_{3}\left(\stackrel{c}{x}^{3}\right)=5.5341
		\end{array}.
		$$
\end{description}

\section{Conclusion}\label{section6}

\par
In this work, we have achieved two objectives: applying the adaptive method to solve a ML(MO)LPP. Secondly, we have enhanced the S.B. Sinha and S. Sinha algorithm \cite{sinha}, by applying the adaptive method in step 10 instead of the simplex method, another advantage of our algorithm is reword the steps (2)-(9) using the the range reduction map $\xi$, which allows us to write those eight steps in one step (Step 8 in our algorithm).
\par
In this way, we use the principle of the constructive adaptive method to propose and demonstrate in detail an algorithm to solve ML(MO)LPP.  We build the range reduction map in order to define the new sub-optimality estimate suitable for the multilevel case which allowed us to nest the algorithm of the adaptive method in the proposed one.
The convergent results of our algorithms are immediately inferred from those of the adaptive method. Finally, we illustrated the algorithm using a numerical example.
 \appendix
\section{Interval reduction function}\label{IRF}
\begin{definition}
Let $a_{1},a_{1},t$ be a real numbers, such that $a_{1}<t<a_{2}$, and define the functions:
\begin{equation}
\begin{array}{cccc}
L:    &     [a_{1},a_{2}]    &  \longrightarrow &    [a_{1},t]\vspace{0.2cm}\\
      &       x      &  \longrightarrow  &    \frac{(t-a_{1})x+a_{1}(a_{2}-t)}{a_{2}-a_{1}}
\end{array},
\end{equation}
\begin{equation}
\begin{array}{cccc}
U:    &     [a_{1},a_{2}]    &  \longrightarrow &    [t,a_{2}]\vspace{0.2cm}\\
      &       x      &  \longrightarrow  &    \frac{(a_{2}-t)x+a_{2}(t-a_{1})}{a_{2}-a_{1}}
\end{array}.
\end{equation}
\end{definition}
\begin{lemma}
The functions $L$ and $U$ are bijectives from $[a_{1},a_{2}]$ to $[a_{1},t]$ $($resp $[t,a_{2}]$ $)$.
\end{lemma}
\begin{proof}
The derivatives of $L^{'}$ and $U^{'}$ are constants and  different from 0 over any interval  $[a_{1},a_{2}]$ of $\mathbb{R}$, and 
$$
\begin{array}{ccc}
L\left([a_{1},a_{2}]\right)=[a_{1},t]&\text{and}&U\left([a_{1},a_{2}]\right)=[t,a_{1}].
\end{array}
$$
\end{proof}
\section{The sign application}\label{IRF1}
\begin{definition}
The application sign is defined as follows:
$$
\begin{array}{ccccc}
\text{sign:}&\mathbb{R}&\longrightarrow&\left\{-1, 0,+1\right\}&\vspace{0.25cm}\\
     & t        &\longrightarrow&\text{sign}(t):=&\left\{\begin{array}{rrr}
		                                 -1&\text{if}& t<0\\
																		0&\text{if}& t=0\\
																		1&\text{if}& t>0
																				      \end{array}\right. .
\end{array}
$$
\end{definition}

\end{document}